\documentclass[12pt]{amsart}
\usepackage{amssymb}

\numberwithin{equation}{section}
\theoremstyle{plain}
\newtheorem{thm}{Theorem}[section]
\newtheorem{prop}[thm]{Proposition}
\newtheorem{cor}[thm]{Corollary}
\newtheorem{lem}{Lemma}[section]

\theoremstyle{definition}
\newtheorem{defn}{Definition}[section]

\theoremstyle{remark}
\newtheorem{rem}{Remark}[section]

\def\R{{\bf R}}
\def\N{{\bf N}}

\def\p#1{{\left({#1}\right)}}
\def\b#1{{\left\{{#1}\right\}}}
\def\br#1{{\left[{#1}\right]}}
\def\n#1{{\left\|{#1}\right\|}}
\def\abs#1{{\left|{#1}\right|}}
\def\jp#1{{\left\langle{#1}\right\rangle}}
\def\Im{\operatorname{Im}}

\def\supp{\operatorname{supp}}

\title{A smoothing property of Schr\"odinger equations
\\
in the critical case}

\author[]{Michael Ruzhansky and Mitsuru Sugimoto}
\address{
  Michael Ruzhansky:
  \endgraf
  Department of Mathematics
  \endgraf
  Imperial College of Science, Technology and Medicine
  \endgraf
  180 Queen's Gate, London SW7 2BZ, UK
  \endgraf
  {\it E-mail address} {\rm ruzh@ic.ac.uk}
  \endgraf
  \medskip
  Mitsuru Sugimoto:
  \endgraf
  Department of Mathematics, Graduate School of Science
  \endgraf
  Osaka University
  \endgraf
  Machikaneyama-cho 1-16, Toyonaka, Osaka 560-0043, Japan
  \endgraf
  {\it E-mail address} {\rm sugimoto@math.wani.osaka-u.jp}
  }

\thanks{This work was completed with the aid of
\lq\lq UK-Japan Joint Project Grant" by
\lq\lq The Royal Society" and
\lq\lq Japan Society for the Promotion of Science".
}

\date{\today}

\begin{document}
\maketitle
\section{Introduction}
Suppose $n\geq2$.
It is well known that the solution $u(t,x)\in C^1\p{\R_t;L^2\p{\R^n_x}}$
to the Schr\"odinger equation
\begin{equation}\label{oeq}
\left\{
\begin{aligned}
\p{i\partial_t+\triangle_x}\,u(t,x)&=0\\
u(0,x)&=\varphi(x)\in L^2\p{\R^n_x}
\end{aligned}
\right.
\end{equation}
has the smoothing property
$\langle D\rangle^{1/2}u\in L^2_{loc}\p{\R_t\times\R^n_x}$
by Sj\"olin \cite{Sj}.
Its global version
\begin{equation}\label{GSP}
\n{\sigma(X,D) u}_{L^2\p{\R_t\times\R^n_x}}
\leq C
\n{\varphi}_{L^2(\R^n)},
\end{equation}
where
\[
\sigma(X,D)=|x|^{\alpha-1}|D|^\alpha, \quad 1-n/2<\alpha<1/2,
\]
was proved by Kato and Yajima \cite{KY} ($n\geq3$)
and Sugimoto \cite{Su1} ($n\geq2$).
Estimate \eqref{GSP} with
\[
\sigma(X,D)=\p{1+|x|^2}^{-s/2}|D|^{1/2}, \quad s>1/2,
\]
was also proved by Ben-Artzi and Klainerman \cite{BK} ($n\geq3$)
and Chihara \cite{Ch} ($n\geq2$), and we can never obtain
\eqref{GSP} in the critical case
\[
\sigma(X,D)=|x|^{-1/2}|D|^{1/2}
\]
 (Watanabe \cite{W}).
But we insist, in this paper, we can attain the critical case if
we assume some structure condition on the operator $\sigma(X,D)$.
This case is crucial in dealing with further applications to
the low regularity solutions to nonlinear problems, where the structure
condition corresponds to the null structure of the nonlinearity.
\par
We explain our strategy by generalizing the equation, so that we
can understand the meaning of the structure well. This will also
emphasize geometric qualities responsible for smoothing in the
critical case. We set
\begin{equation}\label{operator}
L_p=p(D)^2,
\end{equation}
where
$p(\xi)\in C^\infty\p{\R^n\setminus0}$ is a positive function,
and is positively homogeneous of degree 1,
that is, $p(\xi)$ satisfies $p(\lambda \xi)=\lambda p(\xi)$ for $\lambda>0$
and $\xi\in\R^n\setminus0$.
The case $p(\xi)=|\xi|$ corresponds to the usual Laplacian $L_p=-\triangle$.
We assume that the Gaussian curvature of the hypersurface
\begin{equation}\label{hyper}
\Sigma_p=\b{\xi;p(\xi)=1}
\end{equation}
never vanishes.
Under the setting above, we consider the following
generalized Schr\"odinger equation
\begin{equation}\label{eq}
\left\{
\begin{aligned}
\p{i\partial_t-L_p}\,u(t,x)&=0\\
u(0,x)&=\varphi(x)\in L^2\p{\R^n_x}.
\end{aligned}
\right.
\end{equation}
Furthermore, let $\b{(x(t),y(t));\,t\in\R}$
be the classical orbit, that is, the solutions
of the ordinary differential equation
\begin{equation}\label{cl}
\left\{
\begin{aligned}
\dot{x}(t)&=\p{\nabla_\xi p^2}(\xi(t)), \quad\dot{\xi}(t)=0,
\\
x(0)&=0,\quad \xi(0)=k,
\end{aligned}
\right.
\end{equation}
and consider the set of the path of all
classical orbits
\begin{equation}\label{orbits}
\begin{aligned}
\Gamma_p
&=\b{\p{x(t),\xi(t)};\, t\in\R,\, k\in\R^n\setminus 0}
\\
&=\b{\p{\lambda\nabla p(\xi),\xi};\,
 \xi\in\R^n\setminus 0,\,\lambda\in\R}.
\end{aligned}
\end{equation}
For the symbol $\sigma(x,\xi)$ of pseudo-differential operator $\sigma(X,D)$,
we use the notation
\[
\sigma(x,\xi)\sim |x|^a|\xi|^b
\]
if the symbol $\sigma(x,\xi)$ is smooth in $x\neq0$, $\xi\neq0$,
positively homogeneous of order $a$ with respect to $x$,
and of order $b$ with respect to $\xi$.
Then we have the following main result:
\par
\medskip
\begin{thm}\label{main}
Let $n\geq2$.
Suppose $\sigma(x,\xi)\sim|x|^{-1/2}|\xi|^{1/2}$.
Suppose also the structure condition
\begin{equation}\label{str}
 \sigma(x,\xi)=0\quad \text{if}\quad (x,\xi)\in
 \Gamma_p\quad \text{and}\quad x\neq0.
\end{equation}
Then the solution $u$ to $(\ref{eq})$ satisfies estimate \eqref{GSP}, i.e.
$$
\n{\sigma(X,D) u}_{L^2\p{\R_t\times\R^n_x}}
\leq C
\n{\varphi}_{L^2(\R^n)}.
$$
\end{thm}
\medskip
\par
In Section 5 we will discuss the sharpness of the structure condition
\eqref{str}. We will also present results for operators $L_p$ of
arbitrary orders $m\in \N$ as well as results for first and second
order hyperbolic equations.

The proof of Theorem \ref{main} is carried out by replacing
$\sigma(X,D)$ satisfying structure condition by another operator
$\Omega(X,D)$ having good commutability properties (Lemma
\ref{Com}).
This idea can be realized due to a recent progress on
the global $L^2$-boundedness properties of a class of Fourier
integral operators, which was made by authors \cite{RuS}.
With the
aid of it, we prove Proposition \ref{basiclem} which is a key
result in this paper and will be useful in other various
situations.
\par
We remark that we can deduce a global existence result to
derivative nonlinear equations
\begin{equation}
\left\{
\begin{aligned}
\p{i\partial_t-L_p}\,u(t,x)=&\abs{\sigma(X,D)u}^N\\
u(0,x)=\varphi(x),\quad &t\in\R,\, x\in\R^n
\end{aligned}
\right.
\end{equation}
by using Theorem \ref{main} (and its variant).
The structure
condition for the derivative $\sigma(X,D)$ in the nonlinear term
can weaken the regularity assumptions for the initial data
$\varphi$.
We treat this subject in our forthcoming paper
\cite{RS2}.
\par
We should mention here the result of Sugimoto \cite{Su2} which treated the
special case $p(\xi)=|\xi A|$,
where $A$ is a positive definite symmetric matrix, and essentially
proved \eqref{GSP} with
\begin{equation}\label{Su2}
\sigma(X,D)=|x|^{\alpha-1}\Lambda^{1/2-\alpha}
|D|^\alpha,
\quad 1-n/2<\alpha<1/2.
\end{equation}
Here $\Lambda^{\sigma}$ denotes the homogeneous extension of the
(fractional power $\sigma/2$ of) the Laplace-Beltrami operator
on the hypersurface
\begin{align*}
\Sigma_p^*&=\b{\nabla p(\xi);p(\xi)=1}
\\
&=\b{\xi;p^*(\xi)=1}\,;\quad p^*(\xi)=|\xi A^{-1}|.
\end{align*}
Especially $\Lambda^2$ is a partial differential operator of order $2$
with coefficients of homogeneous functions of order $2$,
whose principal symbol is
\begin{equation}\label{tau}
\tau(x,\xi)=
\sum_{i<j}
\p{\frac{p^*(x)}{\abs{\nabla p^*(x)}}\partial_ip^*(x)\xi_j
   -\frac{p^*(x)}{\abs{\nabla p^*(x)}}\partial_jp^*(x)\xi_i}^2
\end{equation}
(see \cite[p.18]{Su2}).
We remark that $\sigma(X,D)$ given by
\eqref{Su2} behaves like $|x|^{-1/2}|D|^{1/2}$ in the sense of the
orders of differentiation and decay, and $\tau(x,\xi)$ is the
typical example which satisfies the structure condition
\eqref{str} (see Remark \ref{remark2}).
Theorem \ref{main} says
that this condition is sufficient for the smoothing property and
the exact form \eqref{Su2} is not necessary.
Even in this sense,
Theorem \ref{main} provides a new aspect to the result of
\cite{Su2}.
\par
We explain the plan of this paper.
In Section 2, we
introduce a class of Fourier integral operators which is the main
tool for the proof of Theorem \ref{main}.
In Section 3, we
investigate the structure of the hypersurface $\Sigma_p$ defined
by \eqref{hyper}, and show the key result Proposition
\ref{basiclem} which is associated to the structure.
In Section 4,
we prove a refined version of the limiting absorption principle which
is an essential part of the proof of Theorem \ref{main}.
In Section 5, we complete the proof and present an extended result of
Theorem \ref{main}.
\par
Finally we remark that the capital $C$'s in estimates
always denote unimportant constants, as usual, which depend only
on the operators, function spaces, and specified suffices.
We remark also that vectors are represented as rows.
\par
\section{Fourier Integral Operator}
We show fundamental properties of a class of Fourier integral operators,
which will be used in the following sections.
To have more flexibility,
we will work in cones in the $(y,\xi)$--space.
Let $\Gamma_y,\Gamma_\xi\subset\R^n$ be open cones.
For the amplitude function 
$a(x,y,\xi)\in C^\infty\p{\R^n_x\times\R^n_y\times\R^n_\xi}$
such that
$\cup_{x\in\R^n}\supp_{y,\xi}a(x,\cdot,\cdot)\subset\Gamma_y\times\Gamma_\xi$
and the phase function
$\varphi(y,\xi)\in C^\infty\p{\Gamma_y\times\Gamma_\xi}$,
we define the operator $T_a$ by
\begin{equation}\label{Ta}
 T_a u(x)
 =\int_{\R^n}\int_{\R^n}
   e^{i\p{x\cdot \xi+\varphi(y,\xi)}}a(x,y,\xi)u(y) dy d\xi.
\end{equation}
In the case $\varphi(y,\xi)=-y\cdot\xi$, $T_a$ is a pseudo-differential
operator, and we use the notation $a(X,Y,D)=\p{2\pi}^{-n}T_a$.
\par
We always assume that $\varphi(y,\xi)$
is a real-valued function and satisfies
\begin{equation}\label{cond:phase}
\begin{aligned}
&\abs{\det\partial_y\partial_{\xi}\varphi(y,\xi)}\geq C>0
\quad\text{on} \quad \Gamma_y\times\Gamma_\xi,
\\
&\abs{\partial_y^\alpha\partial_\xi^\beta\varphi(y,\xi)} \leq
C_{\alpha\beta}\langle y\rangle^{1-|\alpha|}\quad
\text{on}\quad\cup_{x\in\R^n}\supp_{y,\xi}a(x,\cdot,\cdot)
\quad(\forall \alpha,\forall|\beta|\geq1).
\end{aligned}
\end{equation}
Here we have used the notation
\[
\langle x\rangle=\p{1+|x|^2}^{1/2}.
\]
The inequalities
\[
\jp{x}\leq 1+|x|\leq\sqrt{2}\jp{x}
\]
will be used frequently in the following.
From \eqref{cond:phase}, we obtain the estimate
\begin{equation}\label{yequiv}
C_1\jp{y}\leq\jp{\partial_\xi\varphi(y,\xi)}\leq C_2\jp{y}
\quad
\text{on}
\quad
\cup_{x\in\R^n}\supp_{y,\xi}a(x,\cdot,\cdot)
\end{equation}
with some $C_1,C_2>0$.
In fact, the estimate $\jp{\partial_\xi\varphi(y,\xi)}\leq C_2\jp{y}$
is clear.
As for the estimate
$C_1\langle y\rangle\leq \langle \partial_\xi\varphi(y,\xi)\rangle$,
we may assume that $\Gamma_y$ is a proper cone, and fix an appropriate
$y_0\in\Gamma_y$ such that
$\abs{\partial_\xi\varphi(y_0,\xi)}\leq C\jp{y_0}$.
From the expression
\[
\partial_\xi\varphi(y,\xi)-\partial_\xi\varphi(y_0,\xi)
=(y-y_0)\partial_y\partial_\xi\varphi(z,\xi)
\]
with some $z\in\Gamma_y$,
we obtain the estimate
$|y-y_0|\leq C\abs{\partial_\xi\varphi(y,\xi)-\partial_\xi\varphi(y_0,\xi)}$,
hence $\jp{y}\leq C_{y_0}\jp{\partial_\xi\varphi(y,\xi)}$.
\par
For the amplitude function, we introduce classes which emphasize
natural growth properties in variables $x$ and $y$.
\par
\medskip
\begin{defn}\label{class:amp}
Suppose $m,m',k\in\R$.
A smooth amplitude function $a(x,y,\xi)$ is
of the class $\mathcal{A}^{m,m'}_k$, $\mathcal{B}^{m,m'}_k$,
$\mathcal{R}^{m,m'}_k$ respectively if it satisfies
\begin{align*}
&\abs{\partial_x^\alpha\partial_y^\beta\partial_\xi^\gamma
a(x,y,\xi)} \leq C_{\alpha\beta\gamma}
\jp{x}^{m-|\alpha|}\jp{y}^{m'-|\beta|}\jp{\xi}^{k-|\gamma|},
\\
&\abs{\partial_x^\alpha\partial_y^\beta\partial_\xi^\gamma
a(x,y,\xi)} \leq C_{\alpha\beta\gamma}
\jp{x}^{m-|\alpha|}\jp{y}^{m'-|\beta|}\jp{\xi}^{k},
\\
&\abs{\partial_x^\alpha\partial_y^\beta\partial_\xi^\gamma
a(x,y,\xi)} \leq C_{\alpha\beta\gamma}
\jp{x}^{m}\jp{y}^{m'-|\beta|}\jp{\xi}^{k}
\end{align*}
respectively for all $\alpha$, $\beta$, and $\gamma$.
We set $\mathcal{A}^m_k=\cup_{m'\in\R}\mathcal{A}^{m-m',m'}_k$,
$\mathcal{B}^m_k=\cup_{m'\in\R}\mathcal{B}^{m-m',m'}_k$,
$\mathcal{R}^m_k=\cup_{m'\in\R}\mathcal{R}^{m-m',m'}_k$.
In the case $k=0$, we abbreviate it.
\end{defn}
\medskip
\begin{rem}
We have the inclusions
$\mathcal{A}^m_k\subset\mathcal{B}^m_k\subset\mathcal{R}^m_k$.
If the amplitude function is independent of the variable $x$ or $y$,
the definition of these classes can be simplified.
For example, $a(x,\xi)$ is of the class $\mathcal{A}^m$,
$\mathcal{B}^m$ respectively
if it satisfies
\[
\abs{\partial_x^\alpha\partial_\xi^\gamma
a(x,\xi)} \leq C_{\alpha\gamma}
\jp{x}^{m-|\alpha|}\jp{\xi}^{-|\gamma|},
\qquad
\abs{\partial_x^\alpha\partial_\xi^\gamma
a(x,\xi)} \leq C_{\alpha\gamma}
\jp{x}^{m-|\alpha|}
\]
respectively for all $\alpha$ and $\gamma$.
\end{rem}
\medskip
\par
Under the condition \eqref{cond:phase}, we can justify
the definition \eqref{Ta} by the expression
\[
 T_a u(x)
 =\lim_{\varepsilon \searrow 0}
\int_{\R^n}\int_{\R^n}
 e^{i\p{x\cdot \xi+\varphi(y,\xi)}}\rho(\varepsilon\xi)a(x,y,\xi)u(y) dy d\xi
\]
for $a\in\mathcal{R}^{m}_k$,
$u\in C^\infty_0$ and $\rho\in C^\infty_0$.
In fact, by integration by parts argument
\[
e^{i\p{x\cdot \xi+\varphi(y,\xi)}}
=
\frac
 {1-i\p{\partial_y\varphi(y,\xi)-\partial_y\varphi(y,\xi_0)}\cdot\partial_y}
 {1+\abs{\partial_y\varphi(y,\xi)-\partial_y\varphi(y,\xi_0)}^2}
\p{e^{i\p{x\cdot \xi+\varphi(y,\xi)-\varphi(y,\xi_0)}}}e^{i\varphi(y,\xi_0)}
\]
and the inequalities
\begin{align*}
&\abs{\xi-\xi_0}
\leq C\abs{\partial_y\varphi(y,\xi)-\partial_y\varphi(y,\xi_0)},\\
&\abs{\partial_y^\alpha\varphi(y,\xi)-\partial_y^\alpha\varphi(y,\xi_0)}
\leq C_\alpha\jp{y}^{1-|\alpha|}\abs{\xi-\xi_0}
\end{align*}
derived from \eqref{cond:phase} with an appropriate $\xi_0$,
the limit exists and does not depend on the choice of $\rho$.
\par
First we review $L^2$-mapping properties of the operator $T_a$.
For $m\in\R$, let $L^2_m(\R^n)$ be the set of functions $f$ such that the norm
\[
\n{f}_{L^2_m(\R^n)}
=\p{\int_{\R^n}\abs{\langle x\rangle^m f(x)}^2\,dx}^{1/2}
\]
is finite.
Due to Ruzhansky and Sugimoto \cite[Theorem 3.1]{RuS}, we have
\medskip
\begin{thm}\label{Bdd}
Suppose $m\in\R$.
Let $a(x,y,\xi)\in\mathcal{R}^{m}$.
Then
$T_a$ is bounded from $L^2_{m+\mu}(\R^n)$ to $L^2_{\mu}(\R^n)$
for all $\mu\in\R$.
\end{thm}
\medskip
\par
Next we show a symbolic calculus associated to our class.
Before that, we remark the following:
\medskip
\begin{lem}\label{Neg}
Suppose $m,k\in\R$.
For the amplitude function $a(x,y,\xi)\in\mathcal{R}^{m}_k$, we set
\begin{align*}
&a^I(x,y,\xi)=a(x,y,\xi)\chi\p{\p{x+\partial_\xi\varphi(y,\xi)}
     /\langle \partial_\xi\varphi(y,\xi)\rangle},
\\
&a^{II}(x,y,\xi)=a(x,y,\xi)\p{1-\chi}\p{\p{x+\partial_\xi\varphi(y,\xi)}
     /\langle \partial_\xi\varphi(y,\xi)\rangle},
\end{align*}
where $\chi(x)\in C^\infty_0\p{\R^n}$ is equal to one near the origin
and $\supp \chi\subset\b{x;\,|x|<1/2}$.
Then, on $\supp a^I$, we have the equivalence
\[
C_1\jp{x}\leq\jp{y}\leq C_2\jp{x}
\]
for some $C_1,C_2>0$.
Furthermore, for any $N\in\R$, there exists
$r(x,y,\xi)\in\mathcal{R}^{N}_k$
such that $T_{a^{II}}=T_{r}$.
\end{lem}
\medskip
\begin{proof}
Since $\abs{x+\partial_\xi\varphi}\leq(1/2)\jp{\partial_\xi\varphi}$
on $\supp a^I$, we have
\[
\jp{x}
\leq\abs{x+\partial_\xi\varphi}+\sqrt{2}\jp{\partial_\xi\varphi}
\leq (1/2+\sqrt{2})\jp{\partial_\xi\varphi}
\]
and
\begin{align*}
\jp{\partial_\xi\varphi}
&\leq\abs{x+\partial_\xi\varphi}+\sqrt{2}\jp{x}
\\
&\leq (1/2)\jp{\partial_\xi\varphi}+\sqrt{2}\jp{x},
\end{align*}
hence
\[
\jp{\partial_\xi\varphi}\leq2\sqrt{2}\jp{x}.
\]
On account of \eqref{yequiv}, we have the first assertion.
Furthermore, we have
\[
 T_{a^{II}}u(x)
 =\int\int e^{i(x\cdot\xi+\varphi(y,\xi))}
  M^la^{II}(x,y,\xi)u(y) dy d\xi,
\]
where $M$ is the transpose of the operator
\[
{}^tM=\frac{x+\partial_\xi\varphi}{i\abs{x+\partial_\xi\varphi}^2}
    \cdot \partial_\xi
\]
and $l\in\N$.
We have
$\langle \partial_\xi\varphi(y,\xi)\rangle
\leq C\abs{x+\partial_\xi\varphi(y,\xi)}$
on the support of $a^{II}(x,y,\xi)$,
hence we have
\begin{align*}
 &\langle x\rangle
\leq
  \abs{x+\partial_\xi\varphi(y,\xi)}
 +\sqrt{2}\langle \partial_\xi\varphi(y,\xi)\rangle
\leq C\abs{x+\partial_\xi\varphi(y,\xi)},
\\
 &\langle y\rangle
\leq C\langle \partial_\xi\varphi(y,\xi)\rangle
\leq C\abs{x+\partial_\xi\varphi(y,\xi)}
\end{align*}
by \eqref{yequiv}.
On the support of
$\chi'\p{\p{x+\partial_\xi\varphi(y,\xi)}/\jp{\partial_\xi\varphi(y,\xi)}}$,
we have the equivalence of $\abs{x+\partial_\xi\varphi(y,\xi)}$
and $\jp{\partial_\xi\varphi(y,\xi)}$.
By all of these facts, we can justify the second assertion.
\end{proof}
\par
Using this lemma, we have the following fundamental calculus:
\medskip
\begin{thm}\label{Dec}
Suppose $m,k\in\R$.
Let $a(x,y,\xi)\in\mathcal{B}^{m}_k$.
Then we have the decomposition
\[
T_a=T_{a_0}+T_r,
\]
with
\[
a_0(y,\xi)=a(-\partial_\xi\varphi(y,\xi),y,\xi)\in\mathcal{B}^{m}_k,
 \qquad r(x,y,\xi)\in\mathcal{R}^{m-1}_k.
\]
\end{thm}
\medskip
\begin{proof}
The assertion $a_0(y,\xi)\in\mathcal{B}^{m}_k$ is easily obtained
from assumption \eqref{cond:phase}.
We use the decomposition
$a=a^I+a^{II}$ as in Lemma \ref{Neg}, which says that we may
neglect the term $a^{II}$.
By the Taylor expansion, we have
\begin{align*}
a^I(x,y,\xi)
=&a(-\partial_\xi\varphi(y,\xi),y,\xi)
\\
+&\sum_{|\alpha|=1}\p{x+\partial_\xi\varphi(y,\xi)}^\alpha
  \int_0^1\p{\partial_x^\alpha a^I}
   \p{-\partial_\xi\varphi(y,\xi)
+\theta\p{x+\partial_\xi\varphi(y,\xi)},y,\xi}\,
  d\theta.
\end{align*}
Since
\[
\p{x+\partial_\xi\varphi(y,\xi)}^\alpha
e^{i\p{x\cdot\xi+\varphi(y,\xi)}} =\p{-i\partial_\xi}^\alpha
e^{i\p{x\cdot\xi+\varphi(y,\xi)}},
\]
we may take, by integration by parts,
\[
r(x,y,\xi)=
\sum_{|\alpha|=1}
  \int_0^1\p{i\partial_\xi}^\alpha\b{
   \p{\partial_x^\alpha a^I}
   \p{-\partial_\xi\varphi(y,\xi)
+\theta\p{x+\partial_\xi\varphi(y,\xi)},y,\xi}}\
  d\theta.
\]
which belongs to the class $\mathcal{R}^{m-1}_k$
by the equivalence of
$\jp{-\partial_\xi\varphi(y,\xi)+\theta\p{x+\partial_\xi\varphi(y,\xi)}}$
and $\jp{y}$.
\end{proof}
\medskip
\par
By Theorem \ref{Dec} with $\varphi(y,\xi)=-y\cdot\xi$.
we have the following symbolic calculus of pseudo-differential operators:
\medskip
\begin{cor}\label{DecPs}
Suppose $m,k\in\R$.
Let $a(x,y,\xi)\in\mathcal{B}^m_k$.
Then we have the decomposition
\begin{align*}
a(X,Y,D)&=a(Y,Y,D)+r_1(X,Y,D)
\\
&=a(X,X,D)+r_2(X,Y,D),
\end{align*}
with $r_1(x,y,\xi),r_2(x,y,\xi)\in\mathcal{R}^{m-1}_k$.
\end{cor}
\medskip
We now describe a formula for the canonical transform of
pseudo-differential operators.
Let $\Gamma,
\tilde{\Gamma}\subset\R^n\setminus0$ be open cones and
$\psi:\Gamma\to\tilde{\Gamma}$ be a $C^\infty$-diffeomorphism
satisfying $\psi(\lambda\xi)=\lambda\psi(\xi)$ for all
$\xi\in\Gamma$ and $\lambda>0$.
We set formally
\begin{equation}\label{DefI}
\begin{aligned}
&I u(x)
=F^{-1}\left[Fu\p{\psi(\xi)}\right](x)
=(2\pi)^{-n}\int_{\R^n}\int_{\R^n}
   e^{i(x\cdot\xi-y\cdot\psi(\xi))}u(y) dyd\xi,
\\
&I^{-1} u(x)
=F^{-1}\left[Fu\p{\psi^{-1}(\xi)}\right](x)
=(2\pi)^{-n}\int_{\R^n}\int_{\R^n}
  e^{i(x\cdot\xi-y\cdot\psi^{-1}(\xi))}u(y) dyd\xi,
\end{aligned}
\end{equation}
where the Fourier transforms $F$ and $F^{-1}$ are defined by
\[
Fu(x)=\int_{\R^n}e^{-ix\cdot\xi}u(x)\,dx,\qquad
F^{-1}u(\xi)=\p{2\pi}^{-n}\int_{\R^n}e^{ix\cdot\xi}u(\xi)\,d\xi.
\]
The operator $I$ can be justified by using a cut-off function
$\gamma\in C^\infty(\Gamma)$ which satisfies $\supp\gamma\subset\Gamma$,
and $\gamma(\lambda\xi)=\gamma(\xi)$
for large $|\xi|$ and $\lambda\geq1$.
We set
\[
\tilde{\gamma}=\gamma\circ\psi^{-1}\in C^\infty(\tilde{\Gamma}).
\]
Then the operators
\begin{equation}\label{DefI0}
\begin{aligned}
&I_\gamma u(x)
=F^{-1}\left[\gamma(\xi)Fu\p{\psi(\xi)}\right](x)
\\
&=(2\pi)^{-n}\int_{\R^n}\int_{\Gamma}
 e^{i(x\cdot\xi-y\cdot\psi(\xi))}\gamma(\xi)u(y) dyd\xi,
\\
&I_{\tilde{\gamma}}^{-1} u(x)
=F^{-1}\left[\tilde{\gamma}(\xi)Fu\p{\psi^{-1}(\xi)}\right](x)
\\
&=(2\pi)^{-n}\int_{\R^n}\int_{\tilde{\Gamma}}
 e^{i(x\cdot\xi-y\cdot\psi^{-1}(\xi))}\tilde{\gamma}(\xi)u(y) dyd\xi
\end{aligned}
\end{equation}
can be reasonably defined, and we have the expressions
\begin{equation}\label{ipro}
I_\gamma=\gamma(D)\cdot I=I\cdot\tilde{\gamma}(D),\quad
I_{\tilde{\gamma}}^{-1}=\tilde{\gamma}(D)\cdot I^{-1}=I^{-1}\cdot\gamma(D),
\end{equation}
and the identities
\begin{equation}\label{id}
I_\gamma\cdot I_{\tilde{\gamma}}^{-1}=\gamma(D)^2,\quad
I_{\tilde{\gamma}}^{-1}\cdot I_\gamma=\tilde{\gamma}(D)^2.
\end{equation}
We remark that
$\varphi(y,\xi)=-y\cdot\psi(\xi)\in C^\infty(\R^n_y\times\Gamma_\xi)$ and
$\varphi(y,\xi)=
-y\cdot\psi^{-1}(\xi)\in C^\infty(\R^n_y\times\tilde{\Gamma}_\xi)$
satisfy assumption \eqref{cond:phase}.
Then we have the following:
\medskip
\begin{cor}\label{Cantr}
Suppose $m\in\R$.
Let $a(x,\xi)\in\mathcal{A}^m$ (resp. $\mathcal{B}^m$),
and let $I_\gamma$ be defined by
\eqref{DefI0}.
We set $a_0(x,\xi)=a(x,\xi)\gamma_0(\xi)$
with $\gamma_0\in C^\infty(\Gamma)$ satisfying
$\supp\gamma_0\subset\Gamma$ and
$\gamma_0=1$ on $\supp \gamma$, and set
\[
\tilde{a}(x,\xi)=
a_0\p{x \psi'\p{\psi^{-1}(\xi)},\psi^{-1}(\xi)}.
\]
Then we have $\tilde{a}(x,\xi)\in\mathcal{A}^m$ (resp. $\mathcal{B}^m$) and
\[
a(X,D)\cdot I_\gamma
=I_\gamma\cdot\tilde{a}(X,D)+R,
\]
where $R$ is a bounded operator
from $L^2_{m-1+\mu}\p{\R^n}$ to $L^2_\mu\p{\R^n}$ for all $\mu\in\R$.
\end{cor}
\medskip
\begin{proof}
The assertion $\tilde{a}(x,\xi)\in\mathcal{A}^m$ (resp. $\mathcal{B}^m$)
is clear.
Since
\[
b(Y,D)u=F^{-1}_\xi\br{\int e^{-iy\cdot\xi}b(y,\xi)u(y)\,dy},
\]
we have the formula
\[
I_\gamma\cdot b(Y,D)u(x)
=(2\pi)^{-n} \int_{\R^n}\int_{\Gamma} e^{i(x\cdot\xi-y\cdot\psi(\xi))}
  b(y,\psi(\xi))\gamma(\xi)u(y)\,dy d\xi.
\]
Then, by \eqref{ipro} and Theorem \ref{Dec}, we have
\begin{align*}
a(X,D)\cdot I_\gamma u(x)
&=a_0(X,D)\cdot I_\gamma u(x)
\\
&=(2\pi)^{-n} \int_{\R^n}\int_{\Gamma} e^{i(x\cdot\xi-y\cdot\psi(\xi))}
  a_0(x,\xi)\gamma(\xi)u(y)\,dy d\xi
\\
&=(2\pi)^{-n} \int_{\R^n}\int_{\Gamma} e^{i(x\cdot\xi-y\cdot\psi(\xi))}
  a_0\p{y\psi'(\xi),\xi}\gamma(\xi)u(y)\,dy d\xi
+R_1u(x)
\\
&=I_\gamma\cdot \tilde{a}(Y,D)u(x)+R_1u(x),
\end{align*}
and, by Corollary \ref{DecPs}, $\tilde{a}(Y,D)=\tilde{a}(X,D)+R_2$.
We remark that, by Theorem \ref{Bdd}, $I_\gamma$ is bounded
on $L^2_{\mu}\p{\R^n}$, and $R_1$, $R_2$ are bounded
from $L^2_{m-1+\mu}\p{\R^n}$ to $L^2_{\mu}\p{\R^n}$, for all $\mu\in\R$.
We have the corollary by taking $R=R_1+I_\gamma\cdot R_2$.
\end{proof}
\medskip
We now establish the formula for a change of variables.
Let $\Gamma, \tilde{\Gamma}\subset\R^n\setminus0$ be open cones and
$\kappa:\tilde{\Gamma}\to\Gamma$ be a $C^\infty$-diffeomorphism
satisfying $\kappa(\lambda x)=\lambda\kappa(x)$ for all
$x\in\tilde{\Gamma}$ and $\lambda>0$.
We set formally
\begin{equation}\label{DefJ}
\begin{aligned}
&J u(x) =(u\circ\kappa)(x) =(2\pi)^{-n}\int_{\R^n}\int_{\R^n}
e^{i(x\cdot\xi-y\cdot\xi)}
 u\p{\kappa(y)} dyd\xi
\\
&=(2\pi)^{-n}\int_{\R^n}\int_{\R^n} e^{i(x\cdot\xi-\kappa^{-1}(y)\cdot\xi)}
 \abs{\det \p{\kappa^{-1}}'(y)}u(y) dyd\xi,
\\
&J^{-1} u(x) =(u\circ\kappa^{-1})(x)
=(2\pi)^{-n}\int_{\R^n}\int_{\R^n} e^{i(x\cdot\xi-y\cdot\xi)}
 u\p{\kappa^{-1}(y)} dyd\xi
\\
&=(2\pi)^{-n}\int_{\R^n}\int_{\R^n} e^{i(x\cdot\xi-\kappa(y)\cdot\xi)}
 \abs{\det \kappa'(y)}u(y) dyd\xi.
\end{aligned}
\end{equation}
The definition of operator $J$ can be justified by using a cut-off
function $\gamma\in C^\infty(\Gamma)$ which satisfies
$\supp\gamma\subset\Gamma$, and $\gamma(\lambda x)=\gamma(x)$ for
large $|x|$ and $\lambda\geq1$. We set
\[
\tilde{\gamma}=\gamma\circ\kappa\in C^\infty(\tilde{\Gamma}).
\]
Then the operators
\begin{equation}\label{DefJ0}
\begin{aligned}
&J_\gamma u(x)
=\p{\gamma u}\circ\kappa(x)
\\
&=(2\pi)^{-n}\int_{\Gamma}\int_{\R^n} e^{i(x\cdot\xi-\kappa^{-1}(y)\cdot\xi)}
 \abs{\det \p{\kappa^{-1}}'(y)}\gamma(y)u(y) dyd\xi,
\\
&J^{-1}_{\tilde{\gamma}} u(x)
=\p{\tilde{\gamma} u}\circ\kappa^{-1}(x)
\\
&=(2\pi)^{-n}
 \int_{\tilde{\Gamma}}\int_{\R^n} e^{i(x\cdot\xi-\kappa(y)\cdot\xi)}
 \abs{\det \kappa'(y)}\tilde{\gamma}(y)u(y) dyd\xi,
\end{aligned}
\end{equation}
can be reasonably defined, and have the expressions
\begin{equation}\label{jpro}
J_\gamma=J\cdot\gamma(X)=\tilde{\gamma}(X)\cdot J,\quad
J_{\tilde{\gamma}}^{-1}=J^{-1}\cdot\tilde{\gamma}(X)=\gamma(X)\cdot J^{-1}
\end{equation}
and identities
\begin{equation}\label{jd}
J_{\tilde{\gamma}}^{-1}\cdot J_\gamma=\gamma(X)^2,\quad
J_\gamma\cdot J_{\tilde{\gamma}}^{-1}=\tilde{\gamma}(X)^2.
\end{equation}
We remark that 
$\varphi(y,\xi)=-\kappa^{-1}(y)\cdot\xi\in C^\infty(\Gamma_y\times\R^n_\xi)$
and
$\varphi(y,\xi)=-\kappa(y)\cdot\xi\in C^\infty(\tilde{\Gamma}_y\times\R^n_\xi)$
satisfy assumption \eqref{cond:phase}.
Then we have the following:
\medskip
\begin{cor}\label{Chvar}
Suppose $m\in\R$.
Let $a(x,\xi)\in \mathcal{A}^m$, and let $J_\gamma$ be defined
by \eqref{DefJ0}.
We set $a_0(x,\xi)=\gamma_0(x)a(x,\xi)$,
with $\gamma_0\in C^\infty(\Gamma)$
satisfying $\supp\gamma_0\subset\Gamma$ and
$\gamma_0=1$ on $\supp \gamma$, and set
\[
\tilde{a}(x,\xi)=a_0\p{\kappa(x),\xi\kappa'(x)^{-1}}.
\]
Then we have $\tilde{a}(x,\xi)\in \mathcal{A}^m$ and
\[
J_\gamma\cdot a(X,D)
=\tilde{a}(X,D)\cdot J_\gamma+R,
\]
where $R$ is a bounded operator
from $L^2_{m-1+\mu}\p{\R^n}$ to $L^2_{\mu}\p{\R^n}$ for all $\mu\in\R$.
\end{cor}
\medskip
\begin{proof}
The assertion $\tilde{a}(x,\xi)\in \mathcal{A}^m$ is clear.
By \eqref{jpro} and Corollary \ref{DecPs},
we have the decomposition
$J_\gamma\cdot a(X,D)
    =J_\gamma\cdot a_0(X,D)=J_\gamma\cdot a_0(Y,D)+J_\gamma\cdot R_1$.
By Lemma \ref{Neg}, $a_0(y,\xi)$ is also decomposed into the sum
of
\[
a_0^{I}(x,y,\xi)
=a_0(y,\xi)
  \chi\p{\p{x-y}/\jp{y}}
=a_0(y,\xi)
  \chi\p{\p{x-y}/\jp{y}}\gamma_1(y),
\]
where $\gamma_1\in C^\infty(\Gamma)$ satisfies $\supp\gamma_1\subset\Gamma$,
$\gamma_1=1$ on $\supp\gamma_0$,
and $a_0^{II}(x,y,\xi)$ which is negligible.
Furthermore, we have
\begin{align*}
&J_\gamma \cdot a_0^I(X,Y,D)u(x)
\\
&=(2\pi)^{-n}\int\int e^{i(\kappa(x)-\kappa(y))\cdot\xi}
 \gamma(\kappa(x))a_0^I(\kappa(x),\kappa(y),\xi)\abs{\det \kappa'(y)}
 \p{\gamma_1u}\p{\kappa(y)}\,dyd\xi
\\
&=b(X,Y,D)\cdot J_{\gamma_1}u(x),
\end{align*}
where
\begin{align*}
b(x,y,\xi)=&\gamma(\kappa(x))a_0(\kappa(y),\xi\Phi(x,y)^{-1})
              \chi\p{\p{\kappa(x)-\kappa(y)}/\jp{\kappa(y)}}
\\
           &\cdot\abs{\det \kappa'(y)}\cdot\abs{\det\Phi(x,y)}^{-1},
\\
\Phi(x,y)=&\int^1_0\kappa'(x+\theta(y-x))\,d\theta.
\end{align*}
We have used here $\kappa(x)-\kappa(y)=(x-y){}^t\Phi(x,y)$.
Taking $\supp\chi$ to be sufficiently small, we have
\[
|x-y|=\abs{\kappa^{-1}\p{\kappa(x)}-\kappa^{-1}\p{\kappa(y)}}
\leq C\abs{\kappa(x)-\kappa(y)}
\leq \varepsilon\jp{\kappa(y)}
\]
on $\supp b$, where $\varepsilon>0$ is small.
We also have the equivalence of $\jp{\kappa(x)}$ and
$\jp{\kappa(y)}$ by Lemma \ref{Neg}, and
$\jp{\kappa(x)}\leq C\jp{x}$, $\jp{\kappa(y)}\leq C\jp{y}$.
Hence we have the estimates
\[
\jp{x}\leq\ C\jp{x+\theta(y-x)},\quad\jp{y}\leq C\jp{x+\theta(y-x)}
\]
which imply $b(x,y,\xi)\in \mathcal{A}^m$.
Then, by Corollary \ref{DecPs} again, we have $b(X,Y,D)=b(Y,Y,D)+R_2$ and
$b(Y,Y,D)\cdot J_{\gamma_1}=\tilde{a}(Y,D)\cdot J_{\gamma}
=\tilde{a}(X,D)\cdot J_{\gamma}+R_3\cdot J_\gamma$.
By Theorem \ref{Bdd},
we have the corollary with
$R=J_\gamma\cdot R_1 +J_\gamma\cdot a_0^{II}(X,Y,D)+R_2\cdot J_{\gamma_1}
+R_3\cdot J_\gamma$.
\end{proof}
\par
\section{Geometrical structure}
We will now show some basic facts related to the function $p(\xi)$
which is used to define the operator \eqref{operator} in the
introduction.
As we remarked in Introduction, vectors are represented as rows.
We use the notation
\[
\Sigma_q=\b{\xi\in\R^n\setminus0;q(\xi)=1},\quad
\Sigma_q^*=\b{\nabla q(\xi);\xi\in\Sigma_q},
\]
where $\nabla=\partial=\p{\partial_1,\ldots,\partial_n}$.
For vector-valued functions
$f=\p{f_1,\ldots,f_m}$,
we set
$\nabla f={}^t\p{\nabla f_1,\ldots,\nabla f_m}=
\p{\partial_j f_i}_{1\leq i\leq m,1\leq j\leq n}$.
Especially $\nabla^2g=\p{\partial_j\partial_i g}_{1\leq i,j\leq n}$
denotes the Hesse matrix for a scalar valued function $g$.
\medskip
\begin{thm}\label{dual}
Suppose $n\geq2$.
Let $p(\xi)\in C^\infty\p{\R^n\setminus0}$ be a positive and
positively homogeneous function of degree 1.
We assume that the Gaussian curvature of the hypersurface
$\Sigma_p$ never vanishes.
Then there exists a positive and positively homogeneous
function $p^*(\xi)\in C^\infty\p{\R^n\setminus0}$ of order 1
such that
$\Sigma_p^*=\Sigma_{p^*}$,
$\Sigma_{p^*}^*=\Sigma_p$,
and the Gaussian curvature of the hypersurface
$\Sigma_{p^*}$ never vanishes.
Furthermore $\nabla p:\Sigma_p\to\Sigma_{p^*}$ is a $C^\infty$-diffeomorphism
and $\nabla p^*:\Sigma_{p^*}\to\Sigma_p$ is its inverse.
\end{thm}
\medskip
\begin{rem}\label{remark1}
In the case $p(\xi)=\abs{\xi A}$,
where $A$ is a positive definite symmetric matrix,
we have $p^*(\xi)=\abs{\xi A^{-1}}$.
(See \cite[p.19]{Su2}.)
\end{rem}
\medskip
\begin{proof}
First we remark that $\Sigma_p$ is compact.
The curvature
condition is equivalent to the fact that the Gauss map
\[
\Sigma_p\ni\sigma\mapsto\nabla p(\sigma)/\abs{\nabla p(\sigma)}\in S^{n-1},
\]
hence the map $\nabla p:\Sigma_p\to\Sigma_{p^*}$,
is a global $C^\infty$-diffeomorphism.
See Kobayashi and Nomizu \cite{KN},
or consult Matsumura \cite[Theorem D.G. I, p.341]{M}.
Hence, for $\xi\in\R^n\setminus0$,
there is uniquely determined $\sigma\in \Sigma_p$ such that
$\xi/|\xi|=\nabla p(\sigma)/\abs{\nabla p(\sigma)}$, and
we set $\tau(\xi)=\nabla p(\sigma)$.
We remark that $\tau(\xi)=\R_+\xi\cap\Sigma_p^*$,
where $\R_+\xi=\{\lambda\xi;\lambda>0\}$.
Then the positive function $p^*(\xi)=|\xi|/|\tau(\xi)|$ is
in $C^\infty(\R^n\setminus0)$ and
satisfies $p^*(\lambda\xi)=\lambda p^*(\xi)$ for $\lambda>0$
and $\xi\in\R^n\setminus0$.
Since $\xi$ and $\tau(\xi)$ is in the same direction,
$p^*(\xi)=1$ is equivalent to $\xi=\tau(\xi)$, which means
$\xi\in\Sigma^*_p$.
Thus we have obtained the relation
\begin{equation}\label{p^*}
p^*\p{\nabla p(\xi)}=1
\end{equation}
and $\Sigma_p^*=\Sigma_{p^*}$.
\par
Every other statement is implied from the Euler's identity
$p(\xi)=\nabla p(\xi)\cdot \xi$.
In fact, by differentiating it,
we have
$\nabla p(\xi)=\xi\nabla^2 p(\xi)+\nabla p(\xi)$,
hence $\xi\nabla^2 p(\xi)=0$,
which means
\begin{equation}\label{Kerp}
\xi\in \operatorname{Ker}\nabla^2 p(\xi).
\end{equation}
On the other hand, by differentiating the relation $p^*\p{\nabla p(\xi)}=1$,
we have
\[
\nabla p^*\p{\nabla p(\xi)}\in \operatorname{Ker}\nabla^2 p(\xi).
\]
Since the curvature condition is also equivalent to
$\operatorname{rank} \nabla^2 p(\xi)=n-1$
(see Miyachi \cite[Proposition 1]{Miy}),
we have
\begin{equation}\label{dimKerp}
\dim \operatorname{Ker}\nabla^2 p(\xi)=1.
\end{equation}
Hence $\nabla p^*\p{\nabla p(\xi)}$ is parallel to $\xi$ and we can
write, with a scalar $k(\xi)$,
\[
\nabla p^*\p{\nabla p(\xi)}=k(\xi)\xi.
\]
From this representation and
the Euler's identities $p^*(\xi)=\nabla p^*(\xi)\cdot \xi$,
$p(\xi)=\nabla p(\xi)\cdot \xi$, we have
\begin{align*}
p^*\p{\nabla p(\xi)}
&=\nabla p^*\p{\nabla p(\xi)}\cdot\nabla p(\xi)
\\
&=k(\xi)\xi\cdot\nabla p(\xi)
\\
&=k(\xi)p(\xi).
\end{align*}
By \eqref{p^*}, we have $k(\xi)=1/p(\xi)$ and the relation
\begin{equation}\label{inverse}
\nabla p^*\p{\nabla p(\xi)}=\frac{\xi}{p(\xi)}.
\end{equation}
Furthermore, since
the map $\nabla p:\Sigma_p\to\Sigma^*_p=\Sigma_{p^*}$
is onto, it means that the map
$\nabla p^*:\Sigma_{p^*}\to\Sigma_p$
is its inverse.
Hence we have the relation
\begin{equation}\label{inverse'}
\nabla p\p{\nabla p^*(\xi)}=\frac{\xi}{p^*(\xi)}.
\end{equation}
By the Euler's identity again, we have also
\begin{align*}
p\p{\nabla p^*(\xi)}
&=\nabla p\p{\nabla p^*(\xi)}\cdot \nabla p^*(\xi)
\\
&=\frac{\xi}{p^*(\xi)}\cdot \nabla p^*(\xi)=1,
\end{align*}
and the relation $\Sigma_{p^*}^*=\Sigma_p$.
\end{proof}
\medskip
By Theorem \ref{dual}, especially by \eqref{inverse} and \eqref{inverse'},
we have easily the following corollary:
\medskip
\begin{cor}
Let $p(\xi)\in C^\infty\p{\R^n\setminus0}$ be a function used to define
the operator \eqref{operator}, and let $p^*(\xi)\in C^\infty\p{\R^n\setminus0}$
be the function given by Theorem \ref{dual}.
We set, for $\xi\in\R^n\setminus0$,
\begin{equation}\label{phase}
\psi(\xi)=p(\xi)\frac{\nabla p(\xi)}{\abs{\nabla p(\xi)}},
\qquad
\psi^{-1}(\xi)=|\xi|\nabla p^*(\xi).
\end{equation}
Then $\psi:\R^n\setminus0\to\R^n\setminus0$ is a $C^\infty$-diffeomorphism and
$\psi^{-1}:\R^n\setminus0\to\R^n\setminus0$ is its inverse.
\end{cor}
\medskip
For two vectors $a=(a_1,a_2,\ldots,a_n)$ and $b=(b_1,b_2,\ldots,b_n)$,
we define their outer product as
\[
a\wedge b=(a_ib_j-a_jb_i)_{i<j}.
\]
For $\psi(\xi)$ given by \eqref{phase}, we set
\begin{equation}\label{Omega}
\Omega(x,\xi)=x\psi'(\xi)^{-1} \wedge\psi(\xi)
=\p{\Omega_{ij}(x,\xi)}_{i<j}.
\end{equation}
Then we have
\begin{equation}\label{Omega'}
\Omega(x,\xi)=x\nabla^2p^*\p{\nabla p(\xi)}\wedge p(\xi)\nabla p(\xi).
\end{equation}
In fact, differentiating $\psi^{-1}(\xi)=|\xi|\nabla p^*(\xi)$,
we have
\[
\p{\psi^{-1}}'(\xi)
=|\xi|\nabla^2p^*(\xi)+{}^t\nabla p^*(\xi)\frac{\xi}{|\xi|}.
\]
Hence we have, by \eqref{inverse},
\begin{align*}
\psi'(\xi)^{-1}
&=\p{\psi^{-1}}'\p{\psi(\xi)}
\\
&=|\nabla p(\xi)|\nabla^2 p^*\p{\nabla p(\xi)}
+\frac{{}^t\xi}{p(\xi)|\nabla p(\xi)|} \nabla p(\xi),
\end{align*}
and
\[
x\psi'(\xi)^{-1}
=x|\nabla p(\xi)|\nabla^2 p^*\p{\nabla p(\xi)}
+\frac{x\cdot\xi}{p(\xi)|\nabla p(\xi)|}\nabla p(\xi).
\]
Since $\nabla p(\xi)\wedge\nabla p(\xi)=0$, we have \eqref{Omega'}.
\par
The following is another representation of the set $\Gamma_p$ defined by
\eqref{orbits}.
\medskip
\begin{lem}\label{basicrel}
Let $\Omega(x,\xi)$ be given by \eqref{Omega} with $\psi(\xi)$
in \eqref{phase}.
Then we have the relation
\[
\b{(x,\xi)\in\R^n\times\p{\R^n\setminus0}\,;\,\Omega(x,\xi)=0}
=\b{\p{\lambda\nabla p(\xi),\xi};\xi\in\R^n\setminus 0,\lambda\in\R}.
\]
\end{lem}
\medskip
\begin{proof}
First we remark that, by the same reason as those of
\eqref{Kerp} and \eqref{dimKerp},
we have
\begin{equation}\label{Kerp'}
\eta\in \operatorname{Ker}\nabla^2 p^*(\eta)
\end{equation}
and
\begin{equation}\label{dimKerp'}
\dim \operatorname{Ker}\nabla^2 p^*(\eta)=1.
\end{equation}
We take $\eta=\nabla p(\xi)$.
Then $\Omega(x,\xi)=0$ means that
$x\nabla^2 p^*(\eta)$ is parallel to $\eta$ by \eqref{Omega'},
hence it is equivalent to
\[
x\p{\nabla^2 p^*(\eta)}^2=0
\]
by \eqref{Kerp'} and \eqref{dimKerp'}.
Since $\nabla^2 p^*(\eta)$ is diagonalizable,
it is equivalent to
$x\nabla^2 p^*(\eta)=0$ also, which means
\[
x\in \operatorname{Ker}\nabla^2 p^*(\eta).
\]
By \eqref{Kerp'} and \eqref{dimKerp'} again,
it means that $x$ is pararel to $\eta=\nabla p(\xi)$,
hence we have the lemma.
\end{proof}
\medskip
\begin{lem}\label{Com}
Let $\Omega(x,\xi)$ be given by \eqref{Omega} with $\psi(\xi)$
in \eqref{phase}, and let $h\in C^\infty_0\p{\R\setminus0}$.
Then pseudo-differential operators of the form
$\p{h\circ p}(D)$
satisfies
\[
\left[\Omega_{ij},\p{h\circ p}(D)\right]=0.
\]
\end{lem}
\begin{proof}
By the standard symbolic calculus of pseudo-differential
operators, it is sufficient to show
$\nabla_x\Omega_{ij}\cdot\nabla p(\xi)=0$, which is
\[
\eta\nabla_x\p{x\nabla^2p^*(\eta)}=0
\]
with $\eta=\nabla p(\xi)$ by \eqref{Omega'}.
This can be verified by the fact
$\nabla_x\p{x\nabla^2p^*(\eta)}={}^t\nabla^2 p^*(\eta)=\nabla^2 p^*(\eta)$
and \eqref{Kerp'}.
\end{proof}
\medskip
The following estimate is essential for the proof of Theorem \ref{main}
in Introduction.
\medskip
\begin{prop}\label{basiclem}
Suppose $m\in\R$ and $\varepsilon>0$.
Let $a(x,\xi)\in \mathcal{A}^m$.
We assume that $a(x,\xi)=0$ if
$(x,\xi)\in\Gamma_p$ or $|\xi|<\varepsilon$.
Then we have
\begin{equation}\label{basicest}
\n{a(X,D)u}_{L^2(\R^n)}
\leq C
\p{
\sum_{i<j}\n{\Omega_{ij}(X,D)u}_{L^2_{m-1}(\R^n)}
+\n{u}_{L^2_{m-1}(\R^n)}
}.
\end{equation}
Here $\Gamma_p$ is given by \eqref{orbits}
and $\Omega_{ij}(x,\xi)$ by \eqref{Omega}
with $\psi(\xi)$ in \eqref{phase}.
\end{prop}
\medskip
\begin{rem}\label{remark2}
Let $\tau(x,\xi)$ be the symbol \eqref{tau} in Introduction.
Then $\tau(x,\xi)=0$ if $(x,\xi)\in \Gamma_p$ and $x\neq0$.
If fact, for $x\neq0$, $\tau(x,\xi)=0$ is equivalent to
$\abs{\nabla p^*(x)\wedge\xi}^2=0$, which is true on the set $\Gamma_p$
by \eqref{inverse}.
By cutting it off appropriately, we can construct the example which satisfies
the assumption in Proposition \ref{basiclem} 
\end{rem}
\medskip
\begin{proof}
Let $\gamma\in C^\infty\p{\R^n}$ be a cut-off function which satisfies
$\supp \gamma\subset\R^n\setminus0$ and $\gamma(\xi)=1$ for
$|\xi|\geq\varepsilon/2$,
and $I_\gamma$, $I^{-1}_{\tilde{\gamma}}$ be the operators defined by
\eqref{DefI0}
with the phase functions $\psi(\xi)$, $\psi^{-1}(\xi)$
respectively given by \eqref{phase}.
We remark that
the operators $I_\gamma$ and $I^{-1}_{\tilde{\gamma}}$ are bounded on $L^2_\mu$
for any $\mu\in\R$
by Theorem \ref{Bdd},
and we have also
\[
\Gamma_p=\b{(x,\xi)\in\R^n\times\p{\R^n\setminus0}\,;\,\Omega(x,\xi)=0}
\]
by Lemma \ref{basicrel}.
By Corollary \ref{Cantr}, we have
\begin{align*}
&a(X,D)\cdot I_\gamma
=I_\gamma\cdot\tilde{a}(X,D)+R,
\\
&I_{\tilde{\gamma}}^{-1}\cdot\Omega_{ij}(X,D)\cdot \gamma_0(D)
=\tilde{\Omega}_{ij}(X,D)\cdot I_{\tilde{\gamma}}^{-1}+R',
\end{align*}
where
$\gamma_0\in C^\infty\p{\R^n}$ satisfies
$\supp \gamma_0\subset\R^n\setminus0$ and $\gamma_0(\xi)=1$
on $\supp \gamma$,
\[
\tilde{a}(x,\xi)=a\p{x\psi'\p{\psi^{-1}\p{\xi}},\psi^{-1}(\xi)}
\in\mathcal{A}^m,
\qquad
\widetilde{\Omega}(x,\xi)=x\wedge\xi=\p{\widetilde{\Omega}_{ij}(x,\xi)}_{i<j},
\]
and $R$, $R'$ satisfy
\begin{equation}\label{remainder}
\n{Ru}_{L^2\p{\R^n}}\leq C\n{u}_{L^2_{m-1}\p{\R^n}},
\qquad
\n{R'u}_{L^2_{m-1}\p{\R^n}}\leq C\n{u}_{L^2_{m-1}\p{\R^n}}.
\end{equation}
If we notice \eqref{id}, and use the $L^2_{m-1}$-boundedness of
$\left[\Omega_{ij}(X,D),\,\gamma_0(D) \right]$ and
$\gamma_0(D)$
justified by the symbolic calculus of pseudo-differential operators
and Theorem \ref{Bdd},
estimate \eqref{basicest} is reduced to the estimate
\begin{equation}\label{reducedest}
\n{\tilde{a}(X,D)u}_{L^2(\R^n)}
\leq C
\p{
\sum_{i<j}\n{\widetilde{\Omega}_{ij}(X,D)u}_{L^2_{m-1}(\R^n)}
+\n{u}_{L^2_{m-1}(\R^n)}
}.
\end{equation}
We remark that $\tilde{a}(x,\xi)=0$ on the set
\[
\widetilde{\Gamma}_p
=\b{(x,\xi)\in\R^n\times\p{\R^n\setminus0}\,;\,\widetilde{\Omega}(x,\xi)=0}.
\]
\par
Now we prove \eqref{reducedest}.
By covering $\R^n_x$ using sectors and a ball centered at the
origin, we may assume that
$\supp\tilde{a}(\cdot,\xi)\subset\Gamma=\b{x;\,x_n>|x'|}$.
For the justification, here we have used Theorem \ref{Bdd} and the fact
$\widetilde{\Omega}_{ij}(X,D)$ is transformed, by rotation, to a
linear combination of the elements of $\widetilde{\Omega}(X,D)$.
Furthermore, let $\gamma\in C^\infty\p{\R^n}$ be a cut-off
function which satisfies $\supp \gamma\subset\Gamma$ and
$\gamma=1$ on $\supp\tilde{a}(\cdot,\xi)$, and let $J_\gamma$ be
the operator defined by \eqref{DefJ0} with the change of variable
$\kappa:\tilde{\Gamma}\to\Gamma$, where
$\tilde{\Gamma}=\b{x;\,x_n>\sqrt{2}|x'|}$ and
$\kappa(x)=\p{x',\sqrt{x_n^2-|x'|^2}}$.
We have used the notations
$x=(x_1,x_2\ldots,x_n)$, $x'=(x_1\ldots,x_{n-1})$.
By Corollary
\ref{Chvar}, we have
\begin{align*}
&J_\gamma\cdot \tilde{a}(X,D)=b(X,D)\cdot J_\gamma+R
\\
&J_\gamma\cdot \widetilde{\Omega}_{ij}(X,D)
=
\Theta_{ij}(X,D)\cdot J_\gamma+R',
\end{align*}
where
\[
b(x,\xi)=\tilde{a}(\kappa(x),\xi\kappa'(x)^{-1})\in\mathcal{A}^m,
\qquad
\Theta_{ij}(x,\xi)=\widetilde{\Omega}_{ij}\p{\kappa(x),\xi\kappa'(x)^{-1}},
\]
and $R$, $R'$ satisfy \eqref{remainder}.
Here we have noticed that $\Theta_{ij}(X,D)$ is a differential operator.
We remark that we have
\[
\kappa'(x)^{-1}=
\begin{pmatrix}
E_{n-1} & 0
\\
x'/x_n & \sqrt{x_n^2-|x'|^2}/x_n
\end{pmatrix},
\]
where $E_{n-1}$ is the identity matrix of dimension $n-1$.
By easy computation, we have
\begin{align*}
&\Theta_{ij}(x,\xi)=x_i\xi_j-x_j\xi_i
\qquad(i<j<n),
\\
&\Theta_{in}(x,\xi)=-\sqrt{x_n^2-|x'|^2}\xi_i
\qquad(i<n),
\end{align*}
and $b(x,\xi)=0$ on the set $\b{(x,\xi);\,\xi'=0}$.
Hence estimate \eqref{reducedest} is reduced to
\begin{equation}\label{redest2}
\n{b(X,D)\cdot J_\gamma u}_{L^2(\R^n)}
\leq C
\sum_{i<j}\n{\Theta_{ij}(X,D)\cdot J_\gamma u}_{L^2_{m-1}(\R^n)}
\end{equation}
if we notice \eqref{jd}.
\par
Finally, we prove \eqref{redest2}.
Since
\[
\n{D_i\cdot J_\gamma u}_{L^2_m}
\leq C \n{\Theta_{in}(X,D)\cdot J_\gamma u}_{L^2_{m-1}},
\]
it suffices to show
\[
\n{b(X,D)u}_{L^2(\R^n)}
\leq C
\sum_{i=1}^{n-1}\n{D_iu}_{L^2_m(\R^n)}
\]
for $b(x,\xi)\in\mathcal{A}^m$
which vanishes on the set $\b{(x,\xi);\xi'=0}$.
This can be carried out if we notice the Taylor expansion
\[
b(x,\xi',\xi_n)=b(x,0,\xi_n)+\sum_{i=1}^{n-1}r_i(x,\xi)\xi_i,
\]
where
\[
r_i(x,\xi)=\int_0^1\p{\partial_{\xi_i}b}(x,\theta\xi',\xi_n)\in\mathcal{A}^m,
\]
and the boundedness of $r_i(X,D)$ (Theorem \ref{Bdd}).
\end{proof}
\par
\section{Limiting Absorption Principle}
Let $L_p$ be the operator which appears in equation \eqref{eq}.
First we remark that we can define the operator
\[
K_{d,\chi}=
\p{
L_p-d\mp i0
}^{-1}
\chi\p{D},
\]
the weak limit of
$\p{L_p-d\mp i\varepsilon}^{-1}\chi\p{D}$
as $\varepsilon\searrow0$,
where $d\in\R$ and $\chi\in C^\infty_0\p{\R^n\setminus0}$
(see H\"ormander \cite[Section 14.2]{H}).
Taking $\supp \chi$ away from the points $\xi$ such that $p(\xi)^2=d$,
we can also define the operator $\p{L_p-d\mp i0}^{-1}$ since
\[
\lim_{\varepsilon\searrow0}\p{L_p-d\mp i\varepsilon}^{-1}(1-\chi)\p{D}
\to\p{L_p-d}^{-1}(1-\chi)\p{D}.
\]
For the pseudo-differential operator $\sigma(X,D)$, we set
\[
\sigma(X,D)^*=\overline{\sigma(Y,D)}.
\]
Then we have the following refined version of
the limiting absorption principle:
\medskip
\begin{thm}\label{LAP}
Let $n\geq2$, $d\in\R$, and $\chi\in C^\infty_0\p{\R^n\setminus0}$.
Suppose that $\sigma(x,\xi)\in\mathcal{A}^{-1/2}$ and
$\sigma(x,\xi)=0$ on $\Gamma_p$.
Then we have
\[
\n{
\sigma(X,D)
\p{
L_p-d\mp i0
}^{-1}
\chi\p{D}
\sigma(X,D)^*
v}_{L^2(\R^n)}
\le
C_{d,\chi}
\n{v}_{L^2(\R^n)}.
\]
\end{thm}
\medskip
First we prove Theorem \ref{LAP}.
The argument below is the modified version of the proof of
\cite[Theorem 3.1]{Su2}, and may include the repetition of it.
By Proposition \ref{basiclem} and by taking the adjoint, we may show the
$L^2(\R^n)$-boundedness of the operator
\[
\widetilde{K}_{d,\chi}
=\jp{x}^{-3/2}\p{\Omega_{ij}}^k
 K_{d,\chi}
\p{\Omega_{i'j'}^*}^{k'}\jp{x}^{-3/2},
\]
where $\Omega_{ij}$ is given by \eqref{Omega},
$\Omega_{i'j'}^*$ is the adjoint of $\Omega_{i'j'}$,
and $k,k'=0,1$.
Then $\Omega_{ij}$ almost commutes
with $K_{d,\chi}$ in the sense of Lemma \ref{Com}.
On account of it,
we have the expressions,
\begin{align*}
\widetilde{K}_{d,\chi}
=&\sum_{\nu:finite}
  \jp{x}^{-3/2}\p{L_p-d\mp i0}^{-1}\chi_\nu(D)f_\nu(x)\jp{x}^{-3/2}
\\
=&\sum_{\nu:finite}
    \jp{x}^{-3/2}\tilde{f}_\nu(x)\p{L_p-d\mp i0}^{-1}
    \tilde{\chi}_\nu(D)\jp{x}^{-3/2}
\notag
\end{align*}
where $f_\nu,\tilde{f}_\nu$ are polynomials
of order $2$ at most, and
$\chi_\nu,\tilde{\chi}_\nu\in C_0^\infty$ have their support in
that of $\chi$.
Hence we may assume
\begin{equation}\label{T}
\widetilde{K}_{d,\chi}
=\jp{x}^{-3/2}K_{d,\chi}\jp{x}^{1/2},
\qquad
\widetilde{K}_{d,\chi}
=\jp{x}^{1/2}K_{d,\chi}\jp{x}^{-3/2},
\end{equation}
whichever we need.
\par
We may assume, as well, that $\chi(\xi)\in C^\infty_0(\R^n)$
has its support in a sufficiently small
conic neighborhood of $(0,\ldots,0,1)$.
We split the variables in $\R^n$ in the way of
\[
x=(x',x_n),\, x'=(x_1,\ldots,x_{n-1}).
\]
By the integral kernel representation,
we express the operators $\widetilde{K}_{d,\chi}$ and $K_{d,\chi}$ as
\begin{align*}
\widetilde{K}_{d,\chi}v(x)&=\int \widetilde{K}_{d,\chi}(x,y)v(y)\,dy
\\
&=\int\,dy_n\cdot\int \widetilde{K}_{d,\chi}(x,y)v(y)\,dy',
\end{align*}
\begin{align*}
K_{d,\chi} v(x)&=\int K_{d,\chi}(x,y)v(y)\,dy
\\
&=\int\,dy_n\cdot\int K_{d,\chi}(x,y)v(y)\,dy'.
\end{align*}
The following is fundamental in the proof of the limiting absorption principle:
\medskip
\begin{lem}\label{Res}
Let $\chi(\xi)\in C_0^\infty(\R^n)$ have its support in a small
conic neighborhood of $(0,\ldots,0,1)$.
Then we have
\[
   \n{\int K_{d,\chi}(x,y)v(y)\,dy'}_{L^2(\R_{x'}^{n-1})}
  \leq C_{d,\chi}
  \n{v(\cdot,y_n)}_{L^2(\R^{n-1})},
\]
where $C_{d,\chi}$ is independent of $x_n$ and $y_n$.
\end{lem}
\begin{proof}
We follow the argument in the proof of \cite[Lemma 14.2.1]{H}.
Since we have
\[
K_{d,\chi}(x,y)
=\mathcal{F}^{-1}_{\xi}\br{\p{p(\xi)^2-d\mp i0}^{-1}\chi(\xi)}(x-y),
\]
the integral in the left hand side of the estimate is a partial convolution.
We write $\xi=(\xi',\xi_n)$, $\xi'=(\xi_1,\ldots,\xi_{n-1})$.
By virtue of Plancherel's theorem and the inverse formula,
it suffices to show the boundedness of
\[
F^{-1}_{\xi_n}
  \br{\p{p(\xi)^2-d\mp i0}^{-1}\chi(\xi)}(x_n)
\]
with respect to $x_n$ and small $\xi'$.
We may assume $d>0$.
Let the smooth function $\Xi(\xi',d)$ be defined by the relation
$p\p{\xi',\Xi(\xi',d)}^2=d$, which is uniquely determined
for $\xi'$ near $0$.
We have used here Euler's identity
$p(0,\ldots,0,1)=\partial_{\xi_n}p(0,\ldots,0,1)>0$ and the implicit
function theorem.
Then we have
\[
\p{p(\xi)^2-d\mp i\varepsilon}^{-1}\chi(\xi)
=\p{\xi_n-\Xi(\xi',d)\mp i \varepsilon q(\xi)}^{-1}g(\xi),
\]
where
\[
q(\xi)=\frac{\xi_n-\Xi(\xi',d)}{p(\xi)^2-d}\in C^\infty,
\quad
g(\xi)=\chi(\xi)q(\xi)\in C^\infty_0.
\]
Hence we have
\[
\p{p(\xi)^2-d\mp i0}^{-1}\chi(\xi)
=\p{\xi_n-\Xi(\xi',d)\mp i0}^{-1}g(\xi).
\]
The lemma is just a consequence of the formula
\[
F_\rho^{-1}\left[\frac1{\rho\mp i0}\right](r)
=iY(\pm r)
\]
and the fact that $F_{\xi_n}^{-1}g(\xi)$ is
bounded, where $Y$ denotes the Heaviside function.
\end{proof}
\medskip
\par
By (\ref{T}) and Lemma \ref{Res}, we have easily
\[
   \n{\int\jp{x}^{3/2} \widetilde{K}_{d,\chi}(x,y)
     \jp{y}^{-1/2}v(y)\,dy'}_{L^2(\R_{x'}^{n-1})}
  \leq C_{d,\chi}
  {\n{v(\cdot,y_n)}_{L^2(\R^{n-1})}}
\]
and
\[
   \n{\int\jp{x}^{-1/2} \widetilde{K}_{d,\chi}(x,y)
     \jp{y}^{3/2}v(y)\,dy'}_{L^2(\R_{x'}^{n-1})}
  \leq C_{d,\chi}
  {\n{v(\cdot,y_n)}_{L^2(\R^{n-1})}}.
\]
By interpolation, we have
\[
   \n{\int\jp{x}^{1/2\pm\epsilon} \widetilde{K}_{d,\chi}(x,y)
     \jp{y}^{1/2\mp\epsilon}v(y)\,dy'}_{L^2(\R_{x'}^{n-1})}
  \leq C_{d,\chi}
  {\n{v(\cdot,y_n)}_{L^2(\R^{n-1})}},
\]
hence
\[
   \n{\int \widetilde{K}_{d,\chi}(x,y)v(y)\,dy'}_{L^2(\R_{x'}^{n-1})}
  \leq C_{d,\chi}
  \frac{\n{v(\cdot,y_n)}_{L^2(\R^{n-1})}}
   {|x_n|^{1/2\pm\epsilon}|y_n|^{1/2\mp\epsilon}}
\]
for $0<\epsilon\leq1/2$.
Since $|x_n|^{-2\epsilon}\leq2^{2\epsilon}|x_n-y_n|^{-2\epsilon}$
if $|x_n|\geq|y_n|$ and
$|y_n|^{-2\epsilon}\leq2^{2\epsilon}|x_n-y_n|^{-2\epsilon}$
if $|x_n|\leq|y_n|$,
we have
\[
   \n{\int \widetilde{K}_{d,\chi}(x,y)v(y)\,dy'}_{L^2(\R_{x'}^{n-1})}
  \leq C_{d,\chi}
  \frac{\n{v(\cdot,y_n)}_{L^2(\R^{n-1})}}
   {|x_n|^{1/2-\epsilon}|x_n-y_n|^{2\epsilon}|y_n|^{1/2-\epsilon}}.
\]
Then we have
\begin{align*}
\n{\widetilde{K}_{d,\chi}v}_{L^2(\R^n)}
 &\leq
  \n{
     \int\n{\int \widetilde{K}_{d,\chi}(x,y)v(y)\,dy'}_{L^2(\R_{x'}^{n-1})}\,dy_n
    }_{L^2(\R_{x_n})}
\\
 &\leq
  C_{d,\chi}
  \n{
     \int
        \frac{\n{v(\cdot,y_n)}_{L^2(\R^{n-1})}}
         {|x_n|^{1/2-\epsilon}|x_n-y_n|^{2\epsilon}|y_n|^{1/2-\epsilon}}
     \,dy_n
    }_{L^2(\R_{x_n})}
\\
 &\leq
  C_{d,\chi}\n{v}_{L^2(\R^n)},
\end{align*}
where we have used
the following fact (with the case $n=1$)
proved by Hardy-Littlewood \cite{HL}:
\medskip
\begin{lem}\label{SW}
Suppose $\gamma<n/2$, $\delta <n/2$,
$m<n$, and
$\gamma+\delta+m=n$.
Then we have
\[
\p{
\int_{\R^n}
\abs{
\int_{\R^n}
\dfrac{f\p{y}}{\abs{x}^{\gamma}\abs{x-y}^{m}\abs{y}^\delta}
\,dy
}^2\,dx
}^{1/2}
\le C
\p{ \int_{\R^n}\abs{f(x)}^2\,dx
}^{1/2}.
\]
\end{lem}
\medskip
\noindent
(See also \cite[Theorem B]{SW}.)
Thus we have completed the proof of Theorem \ref{LAP}.
\par
As has been already established,
Theorem \ref{main} is a direct consequence of the following
(see Section 5):
\par
\medskip
\begin{cor}\label{RE}
Suppose $\sigma(x,\xi)\sim|x|^{-1/2}|\xi|^{1/2}$.
Suppose also $\sigma(x,\xi)=0$ if $(x,\xi)\in\Gamma_p$ and $x\neq0$.
Then we have
\[
\sup_{d>0}
\n{
\sigma(X,D)
\p{
L_p-d\mp i0
}^{-1}
\sigma(X,D)^*
v}_{L^2(\R^n)}
\le
C\n{v}_{L^2(\R^n)}.
\]
\end{cor}
\medskip
\begin{proof}
By the scaling argument,
we have only to consider the case $d=1$.
We split the estimate into the following two parts:
\begin{equation}\label{high}
\n{
\sigma(X,D)
\p{
L_p
-1\mp i0
}^{-1}
\p{1-\chi\circ p}(D)
\sigma(X,D)^*
v
}_{L^2\p{\R^n}}
\le
C
\n{v}_{L^2\p{\R^n}},
\end{equation}
\begin{equation}\label{low}
\n{
\sigma(X,D)
\p{
L_p
-1\mp i0
}^{-1}
\p{\chi\circ p}(D)
\sigma(X,D)^*
v
}_{L^2\p{\R^n}}
\le
C
\n{v}_{L^2\p{\R^n}},
\end{equation}
where $\chi(t)\in \mathcal{C}^\infty_0\p{\R_+}$ is a function
which is equal to $1$ near $t=1$.
We use the following lemmas:
\par
\medskip
\begin{lem}\label{lemma1}
Suppose $-n/2<\delta<n/2$.
Then we have
\[
\n{|x|^\delta m(D)u}_{L^2(\R^n)}
\le C
\sum_{|\gamma|\leq n}
\sup_{\xi\in\R^n}
\abs{
|\xi|^{|\gamma|}\partial^\gamma m(\xi)
}
\n{|x|^\delta u
}_{L^2(\R^n)}.
\]
\end{lem}
\begin{proof}
See Kurtz and Wheeden \cite[Theorem 3]{KW}.
\end{proof}
\medskip
\begin{lem}\label{lemma2}
Let $a<n/2$, $0<b<a+n/2$.
Suppose $\tau(x,\xi)\sim|x|^{-a}|\xi|^{-b}$.
Then we have
\[
\n{\tau(X,D)u}_{L^2\p{\R^n}}\leq C \n{|x|^{b-a} u}_{L^2\p{\R^n}}.
\]
\end{lem}
\begin{proof}
We have
\[
\tau(X,D)u(x)=\int K(x,x-y)u(y)\,dy,
\]
where
\[
K(x,z)=F_\xi^{-1}[\tau(x,\xi)](z).
\]
Then we have
\[
\abs{K(x,x-y)}\leq C|x|^{-a}|x-y|^{-(n-b)}.
\]
Hence we have
\begin{align*}
\n{\tau(X,D)|x|^{a-b}u}_{L^2} &\leq
C\n{\int\frac{u(y)}{|x|^a|x-y|^{n-b}|y|^{b-a}}\,dy}_{L^2}
\\
&\leq C\n{u}_{L^2}
\end{align*}
where we have used Lemma \ref{SW}.
\end{proof}
\medskip
\par
First we prove estimate \eqref{high}.
Setting
\[
m(\xi)=|\xi|^2\p{p(\xi)^2-1\mp i0}^{-1}\p{1-\chi\circ p}(\xi)
\]
and
\[
\tau(x,\xi)=\sigma(x,\xi)|\xi|^{-1}\sim|x|^{-1/2}|\xi|^{-1/2},
\]
we have
\[
\sigma(X,D)\p{L_p-1\mp i0}^{-1}\p{1-\chi\circ p}(D)\sigma(X,D)^*
=\tau(X,D)m(D)\tau(X,D)^*.
\]
Estimate \eqref{high} is obtained from Lemma \ref{lemma1} with $\delta=0$
and Lemma \ref{lemma2} with $a=b=1/2$.
\par
Next we prove estimate (\ref{low}).
Let $\rho(x)\in C^\infty_0\p{\R^n}$ be equal to 1 near the origin
and $\tilde{\chi}(t)\in C^\infty_0\p{\R_+}$ be equal to 1 on $\supp \chi$.
We set
\[
\sigma_0(x,\xi)=\rho(x)\sigma(x,\xi),\qquad
\sigma_1(x,\xi)=
\p{1-\rho(x)}\sigma(x,\xi)\p{\tilde{\chi}\circ p}(\xi).
\]
Since $\sigma_1(x,\xi)\in\mathcal{A}^{-1/2}$, we have
\begin{equation}\label{sigma11}
\n{
\sigma_1(X,D)
\p{L_p-1\mp i0}^{-1}\p{\chi\circ p}(D)
\sigma_1(X,D)^*
v
}_{L^2\p{\R^n}}
\le
C
\n{v}_{L^2\p{\R^n}}
\end{equation}
by Theorem \ref{LAP}.
On the other hand, since
\[
\sigma_0(x,\xi)
=|x|^{2\varepsilon}\rho(x)
 \tau_0(x,\xi)
 p(\xi)^{1/2+\varepsilon},
\]
where
\[
\tau_0(x,\xi)=|x|^{-2\varepsilon}\sigma(x,\xi)p(\xi)^{-(1/2+\varepsilon)}
\sim|x|^{-(1/2+2\varepsilon)}|\xi|^{-\varepsilon},
\]
we have, by Lemma \ref{lemma2},
\begin{equation}\label{sigma0}
\begin{aligned}
\n{\sigma_0(X,D)u}_{L^2\p{\R^n}}
&\leq C
\n{\tau_0(X,D)p(D)^{1/2+\varepsilon}u}_{L^2\p{\R^n}}
\\
&\leq C
\n{|x|^{-(1/2+\varepsilon)}p(D)^{1/2+\varepsilon}u}_{L^2\p{\R^n}}
\end{aligned}
\end{equation}
with $0<\varepsilon<(n-1)/4$.
By the same argument, we have
\begin{equation}\label{sigma0'}
\begin{aligned}
&\n{\sigma_0(X,D)\Omega_{ij}u}_{L^2\p{\R^n}}
\\
\leq &C
\p{
\n{|x|^{-(1/2+\varepsilon)}p(D)^{3/2+\varepsilon}u}_{L^2\p{\R^n}}
+
\n{|x|^{-(1/2+\varepsilon)}p(D)^{1/2+\varepsilon}u}_{L^2\p{\R^n}}
}
\end{aligned}
\end{equation}
for $\Omega_{ij}$ give by \eqref{Omega} with the same $\varepsilon$.
In fact, since the symbol of $\Omega_{ij}$ is linear in $x$ and
positively homogeneous of order $1$ in $\xi$ by \eqref{Omega'},
we have
\begin{align*}
\sigma_0(X,D)\Omega_{ij}
&=\sum_{\mu:finite}\rho(X)
\sigma_\mu(X,D)+
\sum_{\nu:finite}\rho(X)
\sigma_\nu(X,D)
\\
&=\sum_{\mu:finite}\tilde{\rho}(X)\tilde{\sigma}_\mu(X,D)p(D)+
\sum_{\nu:finite}\rho(X)
\sigma_\nu(X,D)
\end{align*}
by the symbolic calculus, where $\sigma_\mu(x,\xi)\sim|x|^{1/2}|\xi|^{3/2}$,
$\tilde{\rho}(x)=|x|\rho(x)$,
$\tilde{\sigma}_\mu(x,\xi)=|x|^{-1}\sigma(x,\xi)
p(\xi)^{-1}\sim|x|^{-1/2}|\xi|^{1/2}$,
and $\sigma_\nu(x,\xi)\sim|x|^{-1/2}|\xi|^{1/2}$.
Note that the operators $\tilde{\rho}(X)\tilde{\sigma}_\mu(X,D)$ and 
$\rho(X)\sigma_\nu(X,D)$ play the same role in justifying \eqref{sigma0'}
as $\sigma_0(X,D)=\rho(X)\sigma(X,D)$ does in \eqref{sigma0}.
By estimate \eqref{sigma0},
the estimate
\begin{equation}\label{sigma00}
\n{
\sigma_0(X,D)
\p{L_p-1\mp i0}^{-1}\p{\chi\circ p}(D)
\sigma_0(X,D)^*
v
}_{L^2\p{\R^n}}
\le
C
\n{v}_{L^2\p{\R^n}}
\end{equation}
is reduced to show the estimate
\begin{equation}\label{prev}
\n{
\abs{x}^{\alpha-1}
\p{
L_p
-1\mp i0
}^{-1}
\p{\chi\circ p}(D)
\abs{x}^{\beta-1}
v
}_{L^2\p{\R^n}}
\le
C
\n{v}_{L^2\p{\R^n}}
\end{equation}
with $\alpha=\beta=1/2-\varepsilon$
for any $\chi\in C^\infty_0(\R_+)$.
By Lemma \ref{lemma1},
it is implied by
\[
\n{
\abs{x}^{\alpha-1}
\abs{D}^{\alpha+\beta}
\p{
L_p
-1\mp i0
}^{-1}
\abs{x}^{\beta-1}
v
}_{L^2\p{\R^n}}
\le
C
\n{v}_{L^2\p{\R^n}}
\]
which was proved by \cite[Theorem 1.2]{SuT}
(see also \cite[Theorem 3.1]{Su1}).
Furthermore,
by Proposition \ref{basiclem} and Lemma \ref{Com},
estimate
\begin{equation}\label{sigma01}
\n{
\sigma_0(X,D)
\p{L_p-1\mp i0}^{-1}\p{\chi\circ p}(D)
\sigma_1(X,D)^*
v
}_{L^2\p{\R^n}}
\leq
C\n{v}_{L^2\p{\R^n}}
\end{equation}
is reduced to show the estimate
\begin{align*}
&
\n{
\sigma_0(X,D)
\p{L_p-1\mp i0}^{-1}\p{\chi\circ p}(D)
\p{\Omega_{ij}}^k\jp{x}^{-3/2}
v
}_{L^2\p{\R^n}}
\\=
&
\n{
\sigma_0(X,D)\p{\Omega_{ij}}^k
\p{L_p-1\mp i0}^{-1}\p{\chi\circ p}(D)\jp{x}^{-3/2}
v
}_{L^2\p{\R^n}}
\\
\leq
&C\n{v}_{L^2\p{\R^n}},
\end{align*}
where $k=0,1$.
By estimate \eqref{sigma0'} and the trivial inequality
$\jp{x}^{-3/2}\leq|x|^{-(1/2+\varepsilon)}$ with $\varepsilon\leq1$,
it is also reduced to estimate \eqref{prev}.
By taking the adjoint, the estimate
\begin{equation}\label{sigma10}
\n{
\sigma_1(X,D)
\p{L_p-1\mp i0}^{-1}\p{\chi\circ p}(D)
\sigma_0(X,D)^*
v
}_{L^2\p{\R^n}}
\le
C
\n{v}_{L^2\p{\R^n}}
\end{equation}
is also obtained.
Summing up estimates \eqref{sigma11}, \eqref{sigma00},
\eqref{sigma01}, and \eqref{sigma10}, we have estimate \eqref{low},
and thus the proof of Corollary \ref{RE} has been completed.
\end{proof}
\par
\section{Concluding remarks}
First we explain why Corollary \ref{RE} implies Theorem \ref{main}.
The same argument was used in \cite{Su1} or references cited there.
Setting $A=\sigma(X,D)$, we have
\[
\sup_{\rho>0}\n{A (L_p-\rho^2-i0)^{-1}A^*f}_{L^2(\R^n)}
\leq C\n{f}_{L^2(\R^n)}
\]
by Corollary \ref{RE}.
If we notice the formula
\begin{alignat*}{1}
\n{\hat{f}}_{L^2(\rho\Sigma_p\,;\,\rho^{n-1}d\omega/|\nabla p|)}^2
&=
\int_{\Sigma_p}\abs{\hat{f}(\rho\omega)}^2\,\frac{\rho^{n-1}d\omega}
                           {\abs{\nabla p(\omega)}}
\\
&=
4(2\pi)^{n-1}\rho
\Im\p{(L_p-\rho^2-i0)^{-1}f,f}
\end{alignat*}
(see H\"ormander \cite[Corollary 14.3.10]{H}),
we have the estimate
\begin{equation}\label{FRes}
\n{\widehat{A^*f}}_{L^2\p{\rho\Sigma_p\,;\,\rho^{n-1}d\omega/|\nabla p|}}
\le
C\sqrt{\rho}\,
\n{f}_{L^2(\R^n)},
\end{equation}
where $\rho>0$, $\rho\Sigma_p=\{\rho\omega;\omega\in\Sigma_p\}$ and
$d\omega$ is the standard surface element of the hypersurface $\Sigma_p$
defined by \eqref{hyper}.
Let
\[
T=e^{-itp(D)^2}:\mathcal{S}(\R^n_x)\to\mathcal{S}'(\R_t\times\R^n_x)
\]
be the solution operator to \eqref{eq}.
Then the formal adjoint
$T^*:\mathcal{S}(\R_t\times\R^n_x)\to\mathcal{S}'(\R^n_x)$
is expressed as
\[
T^*\left[v(t,x)\right]=F^{-1}_{\xi}
  \left[
    \p{F_{t,x}v}\p{-p(\xi)^2,\xi}
  \right].
\]
In fact, by Plancherel's theorem, we have
\begin{align*}
\p{T\left[\varphi(x)\right],v(t,x)}_{L^2(\R_t\times\R^n_x)}
=&(2\pi)^{-n}
  \iint
   e^{-itp(\xi)^2}
   \p{F_x\varphi}(\xi)
   \overline{F_x\left[v(t,x)\right](\xi)}
  \,dtd\xi
\\
=&(2\pi)^{-n}
  \int
   \p{F_x\varphi}(\xi)
   \overline{\p{F_{t,x}v}\p{-p(\xi)^2,\xi}}
  \,d\xi
\\
=&\p{\varphi(x),
     F^{-1}_{\xi}
          \left[
            \p{F_{t,x}v}\p{-p(\xi)^2,\xi}
          \right](x)}_{L^2(R^n_x)}.
\end{align*}
By \eqref{FRes} and Plancherel's theorem again,
we have
\begin{alignat*}{1}
\n{T^*A^*v}^2_{L^2(\R^n)}
&=C\n{\p{F_{t,x}A^*v}(-p(\xi)^2,\xi)}^2_{L^2(\R^n_\xi)}
\\
&=C\int^\infty_0
  \p{\int_{\Sigma_p}
   \left|
     \p{F_{t,x}A^*v}\p{-\rho^2,\rho\omega}
   \right|^2
    \,\frac{\rho^{n-1}d\omega}{|\nabla p(\omega)|}}\,d\rho
\\
&\leq C
 \int^\infty_0\rho
     \n{\p{F_tv}(-\rho^2,x)
     }^2_{L^2(\R^n_x)}
    \,d\rho
\\
&\leq C
 \int^\infty_{-\infty}
     \n{\p{F_tv}(\rho,x)
     }^2_{L^2(\R^n_x)}
    \,d\rho
\\
&=C
    \n{v(t,x)
      }^2_{L^2(\R_t\times\R^n_x)}.
\end{alignat*}
Here we have used the change of variables $\xi\mapsto\rho\omega$
($\rho>0,\omega\in\Sigma_p$) and $\rho^2\mapsto\rho$.
Then, by the duality argument, we have
\[
\n{AT\varphi}_{L^2\p{\R_t\times\R^n_x}}
\le
C\,
\Vert \,\varphi\,\Vert_{L^2\p{\R^n}},
\]
which is the required estimate \eqref{GSP}.
\par
If we follow the same argument for the operator
\[
T=e^{\pm itp(D)}p(D)^{-1/2}:
 \mathcal{S}(\R^n_x)\to\mathcal{S}'(\R_t\times\R^n_x),
\]
we have a result for the second order hyperbolic equations
(see \cite[Section 5]{Su2}).
Let $\Dot{H}^{s}$ be the homogeneous Sobolev space defined by the norm
\[
\n{\varphi}_{\Dot{H}^{s}\p{\R^{n}_x}}
=\n{\abs{D}^{s}\varphi}_{L^2\p{\R^{n}_x}}.
\]
Then we have the following:
\medskip
\begin{thm}\label{hyperbolic}
Let $n\geq2$.
Suppose $\sigma(x,\xi)\sim|x|^{-1/2}|\xi|^{0}$.
Suppose also the structure condition \eqref{str}.
Then the solution $w$ to
the problem
\[
\left\{
\begin{aligned}
\p{\partial_t^2+L_p}w(t,x)&=0\\
w(0,x)&=\varphi\in L^2\p{\R^{n}_x}\\
\partial_tw(0,x)&=\psi\in \Dot{H}^{-1}\p{\R^{n}_x}
\end{aligned}
\right.
\]
satisfies the estimate
\[
  \n{\sigma(X,D)w}_
   {L^2\p{\R_t\times\R^{n}_x}}
 \leq
   C\p{
       \n{\varphi}_{L^2\p{\R^{n}_x}}
      +\n{\psi}_{\Dot{H}^{-1}\p{\R^{n}_x}}
      }.
\]
\end{thm}
\medskip
We can also treat
the operator $L_p$ of other orders.
Instead of \eqref{operator}, we set $L_p=p(D)^m$, where
$p(\xi)$ satisfies the same assumption.
Theorem \ref{hyperbolic} is also obtained from
the following result (where we can take $m=1$).
\par
\medskip
\begin{thm}\label{main'}
Let $n\geq2$ and $m\in\N$.
Suppose $\sigma(x,\xi)\sim|x|^{-1/2}|\xi|^{(m-1)/2}$.
Suppose also the structure condition \eqref{str}.
Then the solution $u$ to $(\ref{eq})$ with $L_p=p(D)^m$ satisfies
estimate \eqref{GSP}.
\end{thm}
\medskip
\par
We omit the proof of Theorem \ref{main'} because
it is just a straight forward modification of the argument
in Sections 4 and 5.
Modification of the result in \cite{SuT} is also needed
where we use the assumption $m\in\N$.
\par
We finish this article by mentioning that the structure condition \eqref{str}
seems to be necessary for estimate \eqref{GSP}.
By Theorem \ref{main},
the symbol
\[
\sigma(x,\xi)=|x|^{-1/2}\abs{\p{x/|x|}\wedge\nabla p(\xi)}^2|\xi|^{1/2}
\sim|x|^{-1/2}|\xi|^{1/2}
\]
is a typical example which satisfies estimate \eqref{GSP}.
Assume that there is another non-negative symbol
$\tau(x,\xi)\sim|x|^{-1/2}|\xi|^{1/2}$ which breaks the structure condition
\eqref{str} but satisfies estimate \eqref{GSP}.
Then we have the estimate
\[
\n{\p{\sigma(X,D)+\tau(X,D)} u}_{L^2\p{\R_t\times\R^n_x}}
\leq C
\n{\varphi}_{L^2(\R^n)}.
\]
We remark that $\sigma(x,\xi)$ is non-negative and vanishes
only on the set $\Gamma_p$.
Since $\tau(x,\xi)$ is also non-negative and never
vanishes on $\Gamma_p$ by the homogeneity,
we have the ellipticity of the symbol
$\sigma(x,\xi)+\tau(x,\xi)\sim|x|^{-1/2}|\xi|^{1/2}$.
By constructing the parametrix,
we have the critical estimate
\[
\n{|x|^{-1/2}|D|^{1/2}u}_{L^2\p{\R_t\times\R^n_x}}
\leq C
\n{\varphi}_{L^2(\R^n)}
\]
which is not true, at least, for the ordinary
Schr\"odinger equation (Watanabe \cite{W}).
By justifying this argument, we can expect the conclusion
that any non-negative symbol
$\sigma(x,\xi)\sim|x|^{-1/2}|\xi|^{1/2}$ satisfying estimate \eqref{GSP}
must have characteristic points contained in $\Gamma_p$.
\par

\end{document}